\documentclass[11pt]{preprint}
\usepackage[dvips]{graphicx}
\usepackage{amsmath,amsfonts,amssymb,latexsym,epsfig}
\usepackage{mathrsfs}
\usepackage{verbatim}
\usepackage{latexsym}
\usepackage{amssymb}
\usepackage{graphics}
\usepackage{amsbsy}
\usepackage{enumerate}
\usepackage{times}
\usepackage{amsmath}
\usepackage{amsthm}
\usepackage{amsfonts}
\usepackage{mhequ}
\usepackage{ifthen}
\usepackage{times}
\newtheorem{lemma}{Lemma}[section]
\newtheorem{remark}[lemma]{Remark}
\newtheorem{theorem}[lemma]{Theorem}
\newtheorem{proposition}[lemma]{Proposition}
\newtheorem{corollary}[lemma]{Corollary}
\newtheorem{prop}[lemma]{Proposition}
\newtheorem{definition}[lemma]{Definition}
\newtheorem{assumption}[lemma]{Assumption}
\numberwithin{equation}{section}
\newcommand{\ie}{{\textit{i.e.}\ }}

\newcommand{\nnormt}[2][] {\E \bigl( \sup_{t \in [0,\tau^*]} \| #2 \|_{\alpha}^{#1} \bigr)}
\newcommand{\nnorm}[2][] {\E \bigl( \sup_{t \in [0,T]} \| #2 \|_{\alpha}^{#1} \bigr)}

\newcommand{\HH}{{\mathcal H}}
\newcommand{\KK}{{\mathcal K}}

\newcommand{\R}{{\mathbb R  }}

\def\eps{\varepsilon}
\let\epsilon\varepsilon
\def\eref#1{(\ref{#1})}

\def\x(#1){x^{\eps}_{#1}}
\def\y(#1){y^{\eps}_{#1}}
\def\xe(#1){\eps x^{\eps}_{#1}}

\def\CH{{\cal H}}

\def\cN{{\cal N}}
\def\cO{{\cal O}}
\newcommand\scal[2][ ]{\ifthenelse{\equal{#1}{ }}{\langle#2\rangle}{}
        \ifthenelse{\equal{#1}{b}}{\bigl\langle#2\bigr\rangle}{}
        \ifthenelse{\equal{#1}{B}}{\Bigl\langle#2\Bigr\rangle}{}
        \ifthenelse{\equal{#1}{bb}}{\biggl\langle#2\biggr\rangle}{}
        \ifthenelse{\equal{#1}{BB}}{\Biggl\langle#2\Biggr\rangle}{}}

\let\E\expect

\let\P\prob

\def\L^#1{\mathrm{L}^{\!#1}}

\newcommand{\rme}{{\mathrm{e}}}

\newcommand{\One}{{\mathbf 1}}
\newcommand{\CO}{{\mathcal O}}

\newcommand{\cF}{{\mathcal F}}
\def\tr{\mathop{\mathrm{tr}}}
\begin{document}
\title{Multiscale Analysis for SPDEs with Quadratic Nonlinearities}
\author{D. Bl\"{o}mker, M. Hairer, G.A. Pavliotis}
\institute{Institut f\"ur Mathematik, Universit\"at Augsburg, Germany \and Mathematics Insitute, The
University of Warwick, UK \and  Department of Mathematics, Imperial College London, UK
 \\ \email{dirk.bloemker@math.uni-augsburg.de}\\
\email{m.hairer@warwick.ac.uk} \\
\email{g.pavliotis@imperial.ac.uk}}
 \date{\today}
\maketitle

\begin{abstract}
In this article we derive rigorously amplitude equations for stochastic PDEs with quadratic
nonlinearities, under the assumption that the noise acts only on the stable modes and 
for an appropriate scaling between the distance from bifurcation and the strength of the 
noise. We show that, due to the presence of two distinct timescales in
our system, the noise (which acts only on the fast modes) gets transmitted to
the slow modes and, as a result, the amplitude equation contains both additive and
multiplicative noise.

As an application we study the case of the one dimensional Burgers equation
forced by additive noise in the orthogonal subspace to its dominant modes. The theory
developed in the present article thus allows to explain theoretically some recent
numerical observations from \cite{Ro:03}.
\end{abstract}
%
%
%
%
%
%
\section{Introduction}
\label{sec:intro}

Stochastic Partial differential equations (SPDEs) with quadratic nonlinearities arise in
various applications in physics. As examples we mention the use of the stochastic Burgers
equation in the study of closure models for hydrodynamic turbulence \cite{ChekYakh95} and
the use of the stochastic Kuramoto-Sivashinsky equation 
or similar models \cite{Cu-Ba:95, LaurCueMak96,DB-Gu-Ra:02, Ra-Ma-Li-Mo-Ha-Sa:00} for the modelling of surface
phenomena. Very often SPDEs have two widely separated characteristic
timescales and it is desirable to obtain a simplified equation which governs the evolution
of the dominant modes of the SPDE and captures the dynamics
of the infinite dimensional stochastic system at the slow timescale. The purpose of 
this paper is to derive rigorously such an \emph{amplitude equation} for a quite general
class of SPDEs with quadratic nonlinearities and, furthermore, to obtain sharp error
estimates. 

Consider, as a working example of the class of SPDEs that we will consider in this paper, 
the following variation on the Burgers equation
\begin{equation}\label{e:burgers}
\partial_t u =\partial_x^2 u +u \partial_xu +(1+\gamma)u+ \sigma \phi
\end{equation}
subject to external forcing $\sigma \phi$ and to Dirichlet boundary conditions on 
$[0,\pi]$.  Since we are working  very far from the inviscid regime, the solutions to 
this equation in the absence of forcing would decay quickly to $0$, were it not for the 
extra linear instability $(1+\gamma)u$. The constant $1$ appearing in this term is taken 
equal to the Poincar\'e constant for the Dirichlet Laplacian on $[0,\pi]$ and is designed 
to render the first mode $\sin(x)$ linearly neutral. The constant $\gamma$ therefore 
describes the linearised behaviour of that mode. The aim of this article is to study the 
behaviour of solutions to \eref{e:burgers} for small $\gamma$ over (large) timescales of 
order $\gamma^{-1}$.

It is well known \cite{Bl:05, Cr-Ho:93} 
(see also \cite{Schneider2, Schneider1}  for general results on unbounded domains)  
that in the absence of forcing, the solution to
\eref{e:burgers} is of the type
\begin{equ}[e:ansatz]
u(t,x) = \sqrt \gamma a(\gamma t) \sin(x) + \CO(\gamma),
\end{equ}
where the amplitude $a$ solves the deterministic Landau equation
\begin{equ}
\partial_t a= a - {\textstyle\frac1{12}} a^3\;.
\end{equ}
If the forcing $\phi$ is taken to be white in time (actually, any stochastic process with 
sufficiently good mixing properties would also do), then, provided that $\sigma = \cO(\gamma)$,
the solution to \eref{e:burgers} is still of the type \eref{e:ansatz}, but $a$ now solves 
a stochastic version of the Landau equation:
\begin{equ}
\partial_t a= a - {\textstyle\frac1{12}} a^3 + \tilde \sigma \xi(t)\;,
\end{equ}
where $\xi$ is white noise in time and the constant $\tilde \sigma$ is proportional
one the one hand to the ratio $\sigma / \gamma$ and on the other hand to the size of the
projection of $\phi$ onto the `slow' subspace spanned by the mode $\sin(x)$ \cite{Bl:05}. 
In particular,
one gets $\tilde \sigma = 0$ if the projection of $\phi$ onto that subspace vanishes.

This naturally raises the question of the behaviour of solutions to \eref{e:burgers}
when the external forcing acts only on the orthogonal complement of the
`slow' subspace. Roberts \cite{Ro:03} considered for example noise acting only on 
the second mode $\sin(2x)$. Using  formal expansions relying on centre
manifold type arguments, he derived a reduced model describing the amplitude of the
dominant mode. Moreover he demonstrated numerically that additive noise is capable of
stabilising the dominant mode, \ie the noise eliminates a small linear instability.
In turns out that, in order to have a non-trivial effect on the limiting amplitude equation,
the strength of the noise should be chosen to scale like $\sqrt{\gamma}$, i.e.\ $\sigma = 
\CO(\sqrt\gamma)$. We show that in this case, one has (after integrating against smooth
test functions) 
\begin{equ}[e:ansatz2]
u(t,x) = \sqrt \gamma a(\gamma t) \sin(x) + \CO(\gamma^{5/8})\;.
\end{equ}
To be more precise, we have additional noise terms of order $\sqrt \gamma$
on higher modes that average out when integrated against test functions, \ie
they are small in some appropriate weak (averaged) sense.

The amplitude $a$ solves a stochastic differential equation of Stratonovich type
\begin{equ}[e:ampl2]
d a= (1 +  \delta_1) a\, dt - {\textstyle\frac1{12}} a^3\,dt + \sqrt{\delta_2 + 
\delta_3 a^2}\circ dB(t)\;.
\end{equ}
Here, the constants $\delta_i, \, i=1,2,3$ are proportional to $\sigma^2/\gamma$,
$\sigma^4/\gamma^2$, and $\sigma^2/\gamma$ respectively,
with proportionality constants depending on the exact nature of the noise. The Wiener
process $B$ can be constructed explicitly from the external forcing $\phi$, but
unless $\delta_2 = 0$ it is \textit{not} given by a simple rescaling.

In the particular case, where $\phi(x,t) = \sin(2x)\xi(t)$ with $\xi$ a white noise,
one has $\delta_1 = -\frac{\sigma^2}{88 \gamma}$, $\delta_2 = 0$ and
$\delta_3 =\frac{\sigma^2}{36\gamma}$. Note that $\delta_1$ is negative, so that
if $\sigma^2 > 88\gamma$, the solution to \eref{e:ampl2} converges to $0$ almost surely.
This explains the stabilisation effect observed in \cite{Ro:03}.

In this article, we justify rigorously expressions of the form \eref{e:ansatz2} for PDEs of the
form \eref{e:burgers} and we obtain formulas for the coefficients in the amplitude equation \eref{e:ampl2}. Unlike 
\cite{Bl:05} we are interested in the situation where the noise \emph{does not} act on 
the slow degrees of freedom directly but gets transmitted to them through the nonlinear 
interaction with the fast degrees of freedom. From a technical point of view, one of the 
main novelties of this article is that it provides explicit error bounds on the difference 
between the solution of the original SPDE and the solution of the approximating 
amplitude equation; this is a key requirement in tackling the infinite dimensional 
problem. Thus, our result is stronger in that respect than weak convergence type results 
in the spirit of e.g.\ \cite{EthKur86, Kur73}. Furthermore,
we provide an explicit coupling between the two solutions, which is not trivial in 
the sense that, unlike in the case where the noise acts on the slow variables directly,
the white noise driving the resulting amplitude equation is not a simple rescaling of the 
noise driving the original equation.

Finite dimensional SDEs with quadratic nonlinearities and two characteristic, widely
separated, timescales were analysed systematically by Majda, Timofeyev and Vanden Eijnden
in a series of papers \cite{MTV01,MTV99}. The SDEs that were studied
by these authors can be thought of as finite dimensional approximations of stochastic 
PDEs with quadratic nonlinearities of the form \eqref{e:burgers} (in fact, the authors
consider finite dimensional approximations of deterministic PDEs and they introduce
stochastic effects by replacing the quadratic self-interaction terms of the unresolved
variables by an appropriate stochastic process). In these papers, techniques from the
theory of singular perturbation theory for Markov processes were used to derive
stochastic amplitude equations with additive and/or multiplicative noise, which can be
either stable or unstable. The results obtained by formal multiscale asymptotics can be
in principle justified rigorously using the theorem of Kurtz \cite{Kur73}, see
also \cite[Thm 3.1, Ch. 12]{EthKur86}. However, since these results lack explicit error 
estimates, it is not clear \textit{a priori} whether they can be applied to the 
infinite dimensional situation  that we study in this paper.

The rest of the paper is organised as follows. In Section~\ref{sec:result} we state the
assumptions that we make, and present our main result. Sections~\ref{sec:firstreduct}
to~\ref{sec:fin-red} are devoted to the proof of our main theorem. In Section \ref{sec:burgers} 
we
apply our theory to the stochastic Burgers equation.
%
\section{Notations, Assumptions and Main Result}\label{sec:result}
The main object of study of the present article is the following SPDE written in the form
(cf. \cite{DapZab92}):
\begin{equation}\label{e:main}
d u = \bigl(-Lu +B(u,u) + \nu \eps^2 u\bigr)\,dt + \eps Q\,dW(t)\;.
\end{equation}
Throughout this paper, we make the following assumptions.

\begin{assumption}\label{ass:1}
$L$ is a nonnegative definite self-adjoint operator with compact resolvent
in some real Hilbert space $\HH$.
\end{assumption}

Let $\|\cdot\|$ be the norm and $\scal{\cdot,\cdot}$ be the inner product in $\HH$.
We denote by $\{e_k \}_{k =1}^{\infty}$ and $\{\lambda_k \}
_{k =1}^{\infty}$ an orthonormal basis of eigenfunctions and the corresponding (ordered) 
eigenvalues. We will furthermore assume that

\begin{assumption}\label{ass:2}
The kernel of $L$ is one dimensional, \ie $\lambda_1 =0$ and $\lambda_2 > 0$. 
\end{assumption}

We will use the
notation $\cN:= \mathrm{span}\{ e_1 \}$, and $P_c$ for the orthogonal projection
$P_c:\HH \rightarrow \cN$. Furthermore, $P_s:= I - P_c$. Before stating our assumptions
on the nonlinearity $B$, we introduce the following interpolation spaces. For $\alpha > 0$,
we will denote by $\HH^\alpha$ the domain of $L^{\alpha / 2}$ endowed with the scalar 
product $\scal{u,v}_\alpha = \scal{u,(1+L)^\alpha v}$ and the corresponding norm 
$\|\,\cdot\,\|_\alpha$. Furthermore, we identify $\HH^{-\alpha}$ with the dual of
$\HH^{\alpha}$ with respect to the inner product in $\HH$. With this notation at hand we 
can state our next assumption.
\begin{assumption}\label{ass:3}
There exists $\alpha \in (0,2)$ and $\beta \in (\alpha - 2,\alpha]$ 
such that
$B(u,v) \colon \HH^\alpha \times \HH^\alpha \to \HH^{\beta}$ 
is a bounded symmetric bilinear map with
\begin{equ}
\label{e:assump_B}
P_c B(e_k,e_k) = 0\;,
\end{equ}
for every $k \ge 1$.
\end{assumption}
Finally, we assume that the Wiener process driving equation \eref{e:main} satisfies:
\begin{assumption}\label{ass:4}
$W$ is a cylindrical Wiener process on $\HH$.
The covariance operator $Q$ is symmetric, bounded, commutes  with $L$, and satisfies
\begin{equation}\label{e:projection}
\scal{e_1,Q e_1} = 0\;.
\end{equation}
Furthermore, $Q^2 L^{\alpha - 1}$ is trace class in $\CH$, where the value of $\alpha$ is the same
as in Assumption \ref{ass:3}.
\end{assumption}
\begin{remark}
The scaling in $\eps$ in equations \eref{e:main} and \eref{e:mainres} below is dictated by
the symmetry assumptions \eref{e:assump_B} and \eref{e:projection}. If either of these 
assumptions were to fail, the scaling considered in this article would not yield a meaningful 
limit.
\end{remark}
\begin{remark}
We could easily allow for a deterministic forcing term $\eps f$ acting 
on the fast modes. We omit this for simplicity of presentation.
\end{remark}
We are interested in studying the behaviour of small solutions to \eqref{e:main} on timescales
of order $\eps^{-2}$. To this end, we define $v$ through 
$\eps v(\eps^2 t) = u(t)$, so that $v$ is the solution to
\begin{equation}
\label{e:mainres}
d v = \bigl(-\eps^{-2}Lv +\nu v +\eps^{-1} B(v,v)\bigr)\,dt + \eps^{-1} Q\,dW(t)\;.
\end{equation}
Note that we made an abuse of notation in that the Wiener process $W$ appearing in 
\eref{e:mainres} is actually a rescaled version of the one appearing in \eref{e:main}, but it 
has the same distribution. Now we are ready to state the main 
result of this article.

\begin{theorem}\label{thm:amplitude}
Let $L$, $B$, and $Q$ satisfy Assumptions~\ref{ass:1}--\ref{ass:4}. 
Fix a terminal time $T > 0$, a number $R$, as well as constants $p>0$ and $\kappa > 0$.
Then there exists $C>0$ such that, for every $0< \eps < 1$ and every solution $v$
of (\ref{e:mainres}) with
initial condition $v_0 \in \HH^\alpha$ and $\|v_0\|_\alpha \le R$,
there exists a stopping time $\tau$ and a Wiener process $B$ such that 
$$
\E \sup_{t \in [0,\tau]} \|P_c v(t) - a(t)e_1\|_\alpha^p \le C \eps^{p/4 - \kappa}\;,
\qquad \P(\tau < T) \le C \eps^p\;.
$$
Here, $a(t)$ is the solution to the stochastic amplitude equation
\begin{equation}\label{e:stoch-ampl}
d a(t) = \left( \tilde\nu a(t) - \tilde\eta a(t)^3 \right) dt + \sqrt{\sigma_b + \sigma_a a(t)^2
} \, dB\;,\quad a(0) = \scal{v_0,e_1}\;,
\end{equation}
where the coefficients $\tilde\nu, \, \tilde\eta, \, \sigma_a$ and $\sigma_b$ are given by
equations \eref{e:nu}, \eref{e:eta}, \eref{e:ab}, respectively.

Furthermore, the fast Ornstein-Uhlenbeck process $z(t)$ solving
$$
dz=-\eps^{-2}L z dt+\eps^{-1} Q\,dW(t), \qquad z(0)=P_s v_0\;,
$$
satisfies $\E \sup_{t \in [0,\tau]} \|P_sv(t)-z(t)\|^p_\alpha \le C \eps^{p - \kappa}$.
\end{theorem}
\begin{remark}
In a weak norm in time (for example $H^{-1}$)  one can show that $z(t)$ is well approximated 
by white in time and colored in space noise of order $\epsilon$. Formally, we can write 
\begin{eqnarray*}
z(t) & = & e^{-tL\eps^{-2}} P_s v_0+ \eps^{-1}\int_0^t e^{-(t-s)L\eps^{-2}} QdW(s) 
\\ & \approx & \eps  L^{-1}Q \partial_t W\;.
\end{eqnarray*}
Of course, for small transient timescales of 
order $\cO(\eps^2)$ the initial value $P_sv_0$ of $z(t)$ has a contribution of order 
$\cO(1)$. Thus estimates of the error uniformly in time are out of reach.
\end{remark}
\begin{remark}
A immediate corollary of our result is that, under the assumptions of Theorem
\ref{thm:amplitude}, we can write
$$
\E \sup_{t \in [0, \tau \eps^{-2}]}\| u(t) - \eps a(\eps^2 t) e_1 - \eps R(t) \|^p_{\alpha}
\leq C \eps^{\frac{5 p}{4} - \kappa},
$$
where $u(t)$ is the solution to \eqref{e:main} with $u(0) = \cO(\eps)$, $a(t)$ is the
solution to the amplitude equation \eqref{e:stoch-ampl} with $\eps a(0) = \langle u(0), e_1
\rangle$ and $R(t) = z (\eps^2 t)$ is the solution to
$$
d R = - L R + Q d W, \quad \eps R(0) = P_s u(0).
$$
The noise that appears in the equation for $R$ is a rescaled version of the noise that
appears in the equation for $z$.
\end{remark}
Let us discuss briefly the main steps in the proof of this result. We first decompose the 
solution of \eqref{e:main} into a slow and a fast part:
\begin{equ}[e:decomposition]
v(t) = P_c v(t) + P_s v(t) =: x(t) + y(t)\;,
\end{equ}
to obtain a system of SDEs for $(x, y)$, equation \eqref{e:x-v}. Our next step is to
apply It\^{o}'s formula to suitably chosen functions of $x$ and $y$ in order to
eliminate the $\cO(1/\eps)$ terms from \eqref{e:x-v}. We furthermore show that we can replace 
the fast process $y$ by an appropriate Ornstein-Uhlenbeck process $z$. In this way, we obtain 
an SDE for $x$ that involves only $x$ and the (infinite-dimensional) Ornstein-Uhlenbeck process $z$. 
This is done in Section \ref{sec:firstreduct}, see Proposition \ref{prop:1st-reduct}.

A general averaging result (with error estimates) for deterministic integrals that involve 
monomials of the infinite dimensional OU process $z$, see Corollary \ref{cor:homog}, enables us 
to eliminate or simplify various terms in the equation for $(x, z)$ and to reduce the evolution 
of $x$ to the integral equation 
\begin{equation}\label{e:1st2ndred}
x(t) = x(0)+ \tilde \nu  \int_0^t x(s)\,ds - \tilde{\eta} \int_0^t \bigl(x(s)\bigr)^3\,ds + M(t) + R(t)\; ,
\end{equation}
where  $R(t)=\cO(\eps^{1/2-\kappa})$ (for arbitrary $\kappa > 0$)
and $M(t)$ is a martingale whose quadratic variation has
an explicit expression in terms of $(x,z)$. (We made an abuse of notation here and wrote
$x^3$ for what should really be $\scal{x,e_1}^3 e_1$.) This is done in 
Section \ref{sec:secred1}.  

The final step in the reduction procedure is to show that the martingale $M(t)$ can be
approximated (pathwise) by the stochastic integral
$$
\tilde{M}(t) = \int_0^t \sqrt{\sigma_b + \sigma_a a^2(s) } \, d B(s)\;,
$$
where $B(t)$ is a suitable one dimensional Brownian motion and $a$ is the solution
to the amplitude equation \eref{e:stoch-ampl}. This is done in 
Sections \ref{susec:approx} and \ref{sec:secred2}. We remark that, whereas the derivation of
equation~\eqref{e:1st2ndred} is independent of the dimensionality of $x(t)$, the third part of 
the proof is valid only in the case where the kernel of $L$ is one dimensional. This is the price we have to 
pay in order to obtain rigorous explicit error estimates on the validity of the amplitude 
equation.

Let us comment briefly on the case  $\dim(\ker(L)) =k > 1$ but finite. The technique employed in 
 the proof of Theorem~\ref{thm:amplitude} would still apply to this problem to give
{\emph weak convergence} of the projection of the solution of \eqref{e:mainres} to 
$\sum_{j=1}^k a_j(t) e_j$, where $a(t)$ would satisfy a vector valued amplitude equation. 
It does not seem possible, however, to obtain pathwise convergence using our approach, since the time 
change employed in the proof of Lemma~\ref{lem:mart} works only in one dimension.  
Neither does it seem straightforward to modify the present proof in such a way that one can
obtain explicit error estimates without using Lemma~\ref{lem:mart}. In the 
case where the amplitude $a(t)$ is a Brownian motion on $\R^k$,  one could consider one 
dimensional projections, as was done in \cite{HairPavl04}. It is not clear however how to 
adapt the argument used in that paper to the case where the amplitude $a(t)$ is the solution of 
a general SDE. Furthermore, the error estimate obtained in \cite{HairPavl04} scales like 
$\eps^{c/k^2}$ for some appropriate small constant $c$. This is cleary not optimal.
%

%
%
\section{The Reduction to Finite Dimensions}
\label{sec:firstreduct}
%

Let us fix a terminal time $T$ and constants $\kappa > 0$ and $p > 0$. Note that
these constants are not necessarily the same as the ones appearing in the statement
of Theorem~\ref{theo:main} above, but can get `worse' in the course of the proof.

Note first that one has
\begin{lemma}\label{lem:existunique}
Under assumptions \ref{ass:1}, \ref{ass:3} and \ref{ass:4}, equation \eref{e:mainres} has a 
unique local (mild) solution $u$ in $\HH^\alpha$ for every $x \in \HH^\alpha$, i.e.\  
$u$ has continuous paths in  $\HH^\alpha$.
\end{lemma}

\begin{proof}
This follows from an application of Picard's iteration scheme
for the mild solution, see for example \cite{DapZab92}. One can check that
the assumption $\beta < \alpha -2$, together with the continuity assumption on $B(\,\cdot\,,\,\cdot\,)$, imply that
the solution map has the required contraction properties for sufficiently small time.
The fact that the stochastic convolution takes values in $\HH^\alpha$ is a consequence
of Assumption~\ref{ass:4}.
\end{proof}

\begin{remark}
Note that we do not rely on a dissipativity assumption of the underlying SPDE \eref{e:mainres}.
Thus we can only establish the existence of local solutions. The existence of solutions 
on a sufficiently long timescale will be shown later to follow from the dissipativity of 
the approximating equations.
\end{remark}

Substituting the decomposition \eqref{e:decomposition} into \eqref{e:mainres}, we
obtain the following system of equations
\minilab{e:x-v}
\begin{equs}
dx &=  \nu x\,dt
    + 2\eps^{-1}P_c B(x,y) \,dt
    + \eps^{-1} P_c B(y,y) \,dt \label{e:gensys1}\\
dy &= \bigl(\nu  - \eps^{-2}L\bigr)y\,dt + \eps^{-1}P_s B(x+y,x+y) \,dt
     + \eps^{-1} Q\,dW(t)\;. \quad \label{e:gensys2}
\end{equs}
Since Lemma~\ref{lem:existunique} does not rule out the possibility of a finite time blow up in 
$\HH^\alpha$ for the
quite general system \eqref{e:x-v}, we introduce the stopping time
\begin{equation}\label{tau*-def}
\tau^* = T \wedge \inf \{t>0 \,|\; \|v(t)\|_\alpha \ge \eps^{-\kappa} \}\;.
\end{equation}
Note that Lemma~\ref{lem:existunique} ensures that, for a fixed initial condition $v_0$
and for $\eps$ sufficiently small, one has $\tau^* > 0$ almost surely.

Let us fix now some notation.
\begin{definition}
For a real-valued family of processes $\{X_\eps(t)\}_{t\ge 0}$ 
we say $X_\eps = \cO(f_\eps)$, if for every $p \ge 1$ there exists a constant $C_p$ such that
\begin{equation}
\label{e:DefO} \E \Big(\sup_{t\in[0,\tau^*)} |X_\eps(t)|^p\Big) \le C_p f_\eps^p \;.
\end{equation}
We say that $X_\eps= \cO(f_\eps)$ (uniformly)  on $[0,T]$, if we can replace the
stopping time $\tau^*$ by the constant time $T$ in \eref{e:DefO}.
If $X_\eps$ is a random variable independent of time, we use the
same notation without supremum in time, i.e.\ $X_\eps= \cO(f_\eps)$ 
if $\E|X_\eps|^p \le C_p f_\eps^p$.

We use the notation $X_\eps = \cO(f_\eps^-)$ if $X_\eps = \CO(f_\eps \eps^{-\kappa})$
for every $\kappa > 0$.
\end{definition}
%

%
\subsection{Approximation of the stable part by an Ornstein-Uhlenbeck process}
%

In this subsection we show that the `fast' process $y(t)$ is actually close to an 
Ornstein-Uhlenbeck process, at least up to time $\tau^*$. We have the following result.
\begin{lemma}\label{lem:approxOU}
Let $z(t)$ be the $\cN^\perp$-valued process solving the SDE
\begin{equation}\label{e:ou}
d z(t) = -\eps^{-2} L z\, dt +\eps^{-1} Q d W(t), \qquad z(0) = y(0)\;.
\end{equation}
Then one has
$\|y( \cdot) - z( \cdot)\|_\alpha = \cO(\eps^{1-})$.
\end{lemma}

\begin{proof}
It follows from the mild formulation of \eref{e:mainres} that
\begin{equation}\label{e:v-ou-diff}
y(t) = z(t) + \frac1\eps \int_0^t
\rme^{-(t-s)L\eps^{-2}}N(x(s),y(s)) \,ds\;,
\end{equation}
where we have used the notation $N(x,y)= \eps\nu y +  P_s B(x+y,x+y)$. From the properties
of $L$ we deduce that there exist positive constants $C, \, c$ such that 
\begin{equ}[e:semigroup]
\bigl\|e^{-Lt}P_s\bigr\|_{\HH^\beta \to \HH^\alpha} \le \left\{\begin{array}{rl} C 
t^{{(\beta 
- \alpha)}/2} & \text{for $t \le 1$,} \\ Ce^{- ct} & \text{for $t \ge 1$.} 
\end{array}\right.
\end{equ}

Since on the other hand Assumption~\ref{ass:3} implies that 
$$\|N(x,y)\|_\beta \le C (1+\|x+y\|_\alpha)^2\;,
$$ 
the claim follows from the
definition of $\tau^*$ and the fact that the right hand side of \eref{e:semigroup} is integrable
for $\beta > \alpha - 2$.
\end{proof}
The above approximation result enables us to obtain estimates on the statistics of the
stopping time $\tau^*$. For this we will need an estimate on the Ornstein-Uhlenbeck
process \eqref{e:ou}.
\begin{lemma} \label{lem:bound-y}
Suppose that Assumption~\ref{ass:4} holds. Then there is a version of $z$ which is 
almost surely $\HH^\alpha$-valued with continuous sample paths.
Furthermore, for every $\kappa_0 > 0$ and every $p > 0$, there exists a constant $C$ such that
\begin{equation}\label{e:ou-bound}
\nnorm[p]{z(t)} \leq C \eps^{- \kappa_0}.
\end{equation}
\end{lemma}
\begin{proof}
It follows for example from the proof of \cite[Theorem~5.9]{DapZab92}.
\end{proof}
An immediate corollary of Lemmas~\ref{lem:approxOU} and~\ref{lem:bound-y} is that
the process $y(t)$ is `almost bounded' in $\epsilon$. More precisely:
\begin{corollary}\label{cor:vbound}
Assume that the conditions of Lemma~\ref{lem:approxOU}  and Lemma~\ref{lem:bound-y}
hold. Then, for every $\kappa_0 > 0$ and every $p > 0$, there exists a constant $C$ such
that
\begin{equation}\label{e:v-bound}
\nnormt[p]{y(t)} \leq C \eps^{- \kappa_0}.
\end{equation}
\end{corollary}
\begin{proof}
Follows from Lemma~\ref{lem:approxOU}, equation \eref{e:ou-bound}, and the triangle inequality.
\end{proof}

Note that the value of $\kappa_0$ appearing in the statement above can be chosen
independently of the value $\kappa$ appearing in the definition of $\tau^*$.
Thus, with high probability, the event $\tau^*<T$ is caused by $x(t)$ getting too large. To be 
more precise:
\begin{corollary}\label{cor:1stBtau*}
Under Assumptions \ref{ass:1} and \ref{ass:3}, for every $p>0$ there exists a constant
$C$ such that
$$\P(\tau^*<T) \le \P\left(| x(\tau^*) |\ge \eps^{-\kappa}\right) + C\eps^p
\quad \hbox{for}\ \eps\in(0,1).
$$
\end{corollary}
\begin{proof}
Follows from Corollary~\ref{cor:vbound} and the Chebyshev inequality.
\end{proof}
%
\subsection{Elimination of the $\cO(\frac1\eps)$ terms}
%

Let us first introduce some notation. Given a Hilbert space
$\HH$ we denote by $\HH \otimes_s \HH$ its symmetric tensor product. 
Similarly, we use the notation $v_1 \otimes_s v_2 = \frac12 \bigl(v_1\otimes v_2 + 
v_2\otimes v_1\bigr)$ for the symmetric tensor product of two elements and 
$(A \otimes_s B)(x\otimes y) = \frac12 \bigl(Ax\otimes By + By\otimes Ax\bigr)$ for the 
symmetric tensor product of two linear operators. 

Let us recall that the scalar product in the tensor product space 
$\HH_\alpha \otimes \HH_\beta$ is given by 
$\scal{u_1 \otimes v_1 , u_2 \otimes v_2}_{\alpha,\beta}
:=\scal{u_1 , u_2}_\alpha\scal{ v_1 , v_2}_\beta $.
With a slight abuse of notation, we  write 
$\scal{\cdot , \cdot}_{\alpha} := \scal{\cdot , \cdot}_{\alpha,\alpha}$.
Furthermore, we extend the bilinear form $B$ to the tensor product space 
by $B( u\otimes v)=B(u,v)$.

With this notation, one can check that\footnote{Since $L$ has a zero eigenvalue, one 
should interpret $(I \otimes_s L)^{-1}$ and $L^{-1}$ as pseudo-inverses,
where, for instance, $L^{-1}=0$ on the kernel $N(L)$.  We will only apply these two
operators to elements in $\cN^\perp$ so that this is of no consequence.}:
\begin{lemma}
\label{lem:3.8}
The operator $(I \otimes_s L)^{-1}$ is
bounded from $\HH_\gamma \otimes_s \HH_\gamma$ to $\HH_{\gamma+1} \otimes_s 
\HH_{\gamma+1}$ for any $\gamma \in \R$.
\end{lemma}

\begin{proof}
It suffices to note that $I \otimes_s L$ is diagonal 
with eigenvalues $(\lambda_j + \lambda_k)/2$ in the basis $e_j\otimes_s e_k$.
Note that $N(I\otimes_s L)=\mathrm{span}(e_1\otimes_s e_1)$.
\end{proof}

Now we are ready to present the main result of this subsection.

\begin{prop}\label{prop:1st-reduct}
Let $x$ and $z$ be as above and let Assumptions~\ref{ass:1}--\ref{ass:4} hold. Then, 
there exists a process $R=  \cO (\eps^{1-})$ such that, for every stopping time $t$ with $t \le \tau^*$ almost
surely, one has
\begin{equs}
x(t) & = x(0) + \nu \int_0^t x\, ds + 4 \int_0^t P_c B\bigl(P_c B(x,z),L^{-1}z\bigr)\,ds \\
&\quad +   2\int_0^t P_c B\bigl(x,L^{-1} P_s B(x+z,x+z)\bigr)\,ds
+ 2\int_0^t P_c B\bigl(P_c B(z,z),L^{-1}z\bigr)\,ds \\
&\quad + \int_0^t  P_c B (I\otimes_s L)^{-1}\bigl(z \otimes_s Q\, dW(t)\bigr)
+ 2\int_0^t P_c B\bigl(x,L^{-1} Q\,dW(s)\bigr) \\
&\quad + \int_0^t  P_c B (I\otimes_s L)^{-1}\bigl(z \otimes_s P_sB(x+z,x+z)\bigr)\, ds + 
R(t)\;. \label{e:ampl1}
\end{equs}
\end{prop}
An immediate corollary is
\begin{corollary}\label{rem:Hoelder}
Under the assumptions of Proposition \ref{prop:1st-reduct}, define the process $x_R(t)$
by $x_R(t) = x(t) - R(t)$. Then for every $p>0$ and every $\tilde\alpha < 1/2$, one has
\begin{equ}
 \sup_{0 \le s < t \le \tau^*} \frac{| x_R(t) - x_R(s)|}{ |t-s|^{\tilde\alpha}} = \cO(\eps^{0-})\;.
\end{equ}
\end{corollary}
\begin{proof}
It follows immediately from \eref{e:ampl1}, using the definition of $\tau^*$. The condition 
$\tilde\alpha<1/2$ arises from the two stochastic integrals in the right hand side of 
\eref{e:ampl1}.
\end{proof}

\begin{proof}[Proof of Proposition \ref{prop:1st-reduct}]
Applying It\^o's formula to $B(x,L^{-1}y)$, we get the following identity in $\CH^\beta$:
\begin{equs}
d B(x,L^{-1}y) &= 2\nu B(x,L^{-1}y)\,dt + 2 \eps^{-1} B\bigl(P_c B(x,y),L^{-1}y\bigr)\,dt \\
& \quad + \eps^{-1} B\bigl(P_c B(y,y),L^{-1}y\bigr)\,dt - \eps^{-2}B(x,y)\,dt \\
&\quad + \eps^{-1}B\bigl(x,L^{-1} P_s B(x+y,x+y)\bigr)\,dt + \eps^{-1} B\bigl(x,L^{-1}
Q\,dW(t)\bigr)\;.
\end{equs}
Combining this with Lemma~\ref{lem:approxOU} and the continuity properties of $B$ stated
in Assumption~\ref{ass:3},
it follows that, for every stopping time $t$ with $t \le \tau^*$ almost surely, one has
\begin{equs}
2\int_0^t B(x,y)\,ds &=  4 \eps \int_0^t B\bigl(P_c B(x,z),L^{-1}z\bigr)\,ds +
2\eps \int_0^t B\bigl(x,L^{-1 }Q\,dW(s)\bigr)\\
&\quad +  2 \eps \int_0^t B\bigl(x,L^{-1} P_s B(x+z,x+z)\bigr)\,ds  \\
&\quad + 2\eps \int_0^t B\bigl(P_c B(z,z),L^{-1}z\bigr)\,ds  + R_1(t)\;, \label{e:Ito1}
\end{equs}
where $R_1(t) = \cO(\eps^{2-})$.

Applying It\^o's formula to $\frac12 (y\otimes y)$, we get the following identity in
$\CH^{\alpha-2}\otimes \CH^{\alpha-2}$:
\begin{equs}
{\textstyle{\frac12}} d(y\otimes y) &= \nu \bigl(y \otimes y\bigr)\, dt - \eps^{-2} \bigl(y 
\otimes_s Ly\bigr)\,dt\\
&\quad + \eps^{-1} y \otimes_s P_sB(x+y,x+y)\, dt + \eps^{-1} y \otimes_s Q\, dW(t)\\
&\quad + \eps^{-2} \sum_{i=1}^\infty Qe_i \otimes Qe_i\, dt\;.\label{e:Ito2}
\end{equs}
Note however that all terms but the second one actually belong to $\CH^{\alpha-1}\otimes \CH^{\alpha-1}$.
Since $B$ is bounded from $\HH^{\alpha} \otimes \HH^\alpha$ into
$\HH^\beta$ and since $(I\otimes_s L)^{-1}$ is bounded from $\HH^{\alpha-1} \otimes 
\HH^{\alpha-1}$
to $\HH^{\alpha} \otimes \HH^\alpha$,
we can aply $P_c B (I\otimes_s L)^{-1}$ to both sides of \eref{e:Ito2}. Noting that
$P_c B(e_i \otimes e_i) = 0$ by Assumption~\ref{ass:3}, we get
\begin{equs}
\int_0^t P_c B(y,y)\,ds &= \eps \int_0^t  P_c B (I\otimes_s L)^{-1}\bigl(z \otimes_s 
P_sB(x+z,x+z)\bigr)\, ds\\
&\quad + \eps \int_0^t  P_c B (I\otimes_s L)^{-1}\bigl(z \otimes_s Q\, dW(t)\bigr) + R_2(t)\;,
\end{equs}
where $R_2(t) = \cO(\eps^{2-})$. Collecting both terms and inserting them into \eref{e:gensys1} 
concludes the proof.
\end{proof}
%
%
%
%
\section{Averaging Over the Fast Ornstein-Uhlenbeck Process}
\label{sec:secred1}
%
%
In this section, we simplify the equation for $x$ further by showing that one can eliminate
all terms in \eref{e:ampl1} that contain odd powers of $z$. Furthermore, concerning the terms that
are quadratic in $z$, there exists a constant $\widehat Q \in \HH^\alpha \otimes_s \HH^\alpha$
so that one can make the formal substitution $z \otimes z \mapsto \widehat Q$.

We start with a number of bounds on the integrals of products of Ornstein-Uhlenbeck processes.
%
%
\subsection{An averaging result with explicit error bounds}
%
Recall that $Q W$ can (at least on a formal level) be written as $\sum_{k=2}^\infty q_k e_k w_k(t)$
for some independent standard Wiener processes $w_k$.
For $\eps > 0$ and $k > 1$, we define $\hat z_k(t)$ to be the stationary solution of
\begin{equ}
d\hat z_k = - \lambda_k \eps^{-2} \hat z_k\,dt + q_k \eps^{-1} dw_k(t)\;.
\end{equ}
This is a Gaussian process with covariance
\begin{equ}[e:covy]
\E \left( \hat z_k(s) \hat z_k(t) \right)  = \frac{q_k^2}{ 2\lambda_k} e^{-\frac{\lambda_k |t-s| }{
\eps^2}}\;.
\end{equ}
Since $\hat z_k(t)$ fluctuates very rapidly, one would expect from the law of large numbers
that as $\eps \to 0$, one has $\hat z_k \to 0$ weakly, but $(\hat z_k)^2 \to \frac{q_k^2}{
2\lambda_k}$ weakly. This is made precise by the following bounds:

\begin{lemma}\label{lem:moments}
For every $p > 0$ there exists a constant $C_p$ such that, for every $t>s> 0$ and every
$k,\ell,m > 1$, the bounds
\begin{equs}
\E \Bigl(\int_s^t \hat z_k(r)\,dr\Bigr)^{2p} &\le C_p \Bigl(\frac{q_k^2}{\lambda_k}\Bigr)^p (t-s)^p\eps^{2p}\;,
\\
\E \Bigl(\int_s^t \Bigl(\hat z_k(r)\hat z_\ell(r)- \frac{q_k^2}{2\lambda_k} \delta_{kl} \Bigr)\,dr\Bigr)^{2p} &\le C_p \Bigl(\frac{q_k^2q_\ell^2}{\lambda_k\lambda_\ell}\Bigr)^p (t-s)^p \eps^{2p}\;,\\
\E \Bigl(\int_s^t \hat z_k(r)\hat z_\ell(r)\hat z_m(r)\,dr\Bigr)^{2p} &\le C_p \Bigl(
\frac{q_k^2q_\ell^2q_m^2}{\lambda_k\lambda_\ell\lambda_m}\Bigr)^p (t-s)^p\eps^{2p}\;,
\end{equs}
hold. Here we denoted by $\delta_{kl}$ the Kronecker symbol.
\end{lemma}

\begin{proof}
The first bound can be checked explicitly in the case $p = 1$ by using \eref{e:covy}.
The case $p > 1$ follows immediately from the fact that $\int_s^t \hat z_k(r)\,dr$ is Gaussian.

In order to obtain the other bounds, we recall first the fact that for an $\R^{2p}$-valued 
Gaussian random variable $X=(X_1,\ldots,X_{2p})$, we have
$$ \E X_1\cdot\ldots\cdot X_n = \sum_{\sigma\in\Sigma(2p)} 
\prod_{i=1}^p \E X_{\sigma_{2i-1}} X_{\sigma_{2i}}
 $$
where  $\Sigma(2p)$ is the set of all permutations
of $\{1,\ldots,2p\}$. 

Turning to the second claim, consider first the case where $k \neq \ell$, so that $\hat z_k$ and $\hat z_\ell$ are independent. Thus
\begin{equs}
\E &\prod_{i=1}^{2p} \hat z_k(t_i)\hat z_\ell(t_i) = \E \prod_{i=1}^{2p} \hat z_k(t_i)
\E \prod_{i=1}^{2p}\hat z_\ell(t_i)\\
& \leq C \Bigl({q_k^2 q_\ell^2 \over \lambda_k\lambda_\ell}\Bigr)^{p} \sum_{\sigma\in\Sigma(2p)}
\prod_{i=1}^p  \exp\Bigl(-\frac{c}{\eps^2}|t_{\sigma_{2i}} - t_{\sigma_{2i-1}}|\Bigr)\;,
\end{equs}
where $c = \lambda_1$ is the smallest non-zero eigenvalue of $L$. The bound then follows by integrating
over $t_1,\ldots,t_{2p}$ and using the fact that $\int_s^t \int_s^t \exp(-c |r-u|/\eps^2)\, dr\,du \le C \eps^2 (t-s)$.

In the case where $k = \ell$, we have
\begin{equs}
\E \prod_{i=1}^{2p} &\Bigl((\hat z_k)^2(t_i) - {q_k^2 \over 2\lambda_k}\Bigr) =
\sum_{A\subset\{1,\ldots,2p\}} \Big(-\frac{q_k^2}{2\lambda_k}\Big)^{2p-|A|}  \E
\prod_{i\in A} \hat z_k(t_i)^2\\
&= \sum_{A\subset\{1,\ldots,2p\}} \Big(-\frac{q_k^2}{2\lambda_k}\Big)^{2p-|A|} 
\sum_{\tau\in\Sigma^2(A)} \prod_{i=1}^{|A|} \E \hat z_k(t_{\tau_{2i}})\hat z_k(t_{\tau_{2i-1}})\;,
\end{equs}
where the sum runs over $\Sigma^2(A)$, the space of all permutations
of numbers in $A$, where each number is allowed to appear twice.
Now it is possible to check that all terms in the double sum where $|A|<2p$ 
are cancelled by a term with a larger $\tilde{A}$, where $t_{\tau_{2i}}=t_{\tau_{2i-1}}$ for some $i$. 
All remaining terms have $|A|=2p$.
It follows from \eref{e:covy} that there exists a constant $C$ such that
\begin{equs}
\E \prod_{i=1}^{2p} \Bigl((\hat z_k)^2(t_i) - {q_k^2 \over 2\lambda_k}\Bigr) & \le C
\Bigl({q_k^2 \over \lambda_k}\Bigr)^{2p} \sum_{\sigma\in\Sigma(2p) \atop \sigma_i \neq i} \exp\Bigl( {-c{|t_{\sigma_1}
- t_{1}| + \ldots  + |t_{\sigma_{2p}} - t_{2p}| \over \eps^{2} }}\Bigr)\;.
\end{equs}
The bound then follows immediately by integrating over $t_1,\ldots,t_p$.
The last term can be bounded in a similar way.
\end{proof}
Let now $\widehat Q \in \HH^\alpha \otimes_s \HH^\alpha$ be given by 
\begin{equ}
\label{e:defQhat}
\widehat Q = \sum_{k =2}^\infty {q_k^2\over 2\lambda_k} \bigl(e_k \otimes e_k\bigr)\;.
\end{equ}
The fact that it does indeed belong to $\HH^\alpha \otimes_s \HH^\alpha$ is a consequence of 
Assumption~\ref{ass:4}. Writing $\hat z(t) = \sum_{k=1}^\infty \hat z_k(t) e_k$, we have the 
following corollary of Lemma~\ref{lem:moments}:
\begin{corollary}\label{cor:averagez}
For every $p > 0$ there exists a constant $C_p$ such that the bounds
\begin{equs}
\E \Bigl\|\int_s^t \hat z(r)\,dr\Bigr\|_\alpha^{2p} &\le C_p (t-s)^p\eps^{2p}\;, \\
\E \Bigl\|\int_s^t \bigl(\hat z(r)\otimes \hat z(r)- \widehat Q\bigr) \,dr\Bigr\|_\alpha^{2p} &\le C_p (t-s)^p \eps^{2p}\;,\\
\E \Bigl\|\int_s^t \bigl(\hat z(r)\otimes \hat z(r)\otimes \hat z(r)\bigr)\,dr\Bigr\|_\alpha^{2p} &\le C_p  (t-s)^p\eps^{2p}\;,
\end{equs}
hold  for every $t>s> 0$.
\end{corollary}

\begin{proof}
It is a straightforward consequence of the following fact. Let $\{v_k\}$ be a sequence of 
real-valued random variables such that there exists a sequence $\{\gamma_k\}$ and, for
every $p\ge1$, a constant $C_p$ such that $\E |v_k|^p \le C_p \gamma_k^p$. Then, 
for every $p \ge 1$ there exists a constant $C_p'$ such that
\begin{equ}
\E \Bigl(\sum_{k =1}^\infty \lambda_k^\alpha v_k^2\Bigr)^p \le C_p'\Bigl(\sum_{k=1}^\infty \lambda^\alpha \gamma_k^2\Bigr)^p\;.
\end{equ}
The result now follows immediately from Lemma~\ref{lem:moments} and from Assumption~\ref{ass:4} which states that
the sequence $q_k^2 \lambda_k^{\alpha-1}$ is summable. 
\end{proof}

\begin{lemma}\label{lem:supav}
Let $G_\eps$ be a family of processes in some Hilbert space $\KK$
such that, for every $p \ge 1$ and every $\kappa > 0$ there exists a constant $C$ such that
\begin{equ}[e:aprioriGen]
\E \Bigl\|\int_s^t G_\eps(r)\,dr\Bigr\|^{2p} \le C (t-s)^p\eps^{2p}\;,
\end{equ}
holds for every $0 \le s < t \le 1$. Then, for every $p > 0$ and every $\kappa > 0$, there exists a
constant $C$ such that
\begin{equ}
\E \Bigl(\sup_{n < N} \Bigl\|\int_{n\delta}^{(n+1)\delta} G_\eps(s)\,ds\Bigr\|\Bigr)^{2p}
\le C N^\kappa \delta^p \eps^{2p}\;,
\end{equ}
holds for every $N>0$, every $\delta \in (0,(N+1)^{-1})$, and every $\eps > 0$.
\end{lemma}

\begin{proof}
It follows from \eref{e:aprioriGen} that, for every $q > 0$, there exists a constant $C$
such that
\begin{equ}[e:boundSum]
\P \Bigl(\sup_{n < N} \Bigl\|\int_{n\delta}^{(n+1)\delta} G_\eps(s)\,ds\Bigr\| > K\Bigr)
\le C N {\delta^{q/2} \over K^q} \eps^{q}\;,
\end{equ}
holds for every $K > 0$. Note now that if a positive random variable $X$ satisfies $\P(X
> x) \le \bar C /x^q$ for every $x > 0$, then, for $p < q$, one has
\begin{equs}
\E X^p &= p\int_0^\infty x^{p-1} \P(X > x)\,dx
\le p\int_0^{\bar C^{1/q}} x^{p-1}\,dx + \bar Cp\int_{\bar C^{1/q}}^\infty x^{p-q-1}\,dx \\
&=  {q \over q-p} \bar C^{p/q}\;.
\end{equs}
Combining this with \eref{e:boundSum} and choosing $q$ sufficiently large yields the
required bound.
\end{proof}

\begin{proposition}\label{prop:homog}
Let $\KK$ be a Hilbert space, let $f$ be a $\KK$-valued random process with almost surely 
$\tilde\alpha$-H\"older continuous trajectories, let $G_\eps$ be a family of $\KK$-valued 
processes satisfying \eref{e:aprioriGen}, and let
\begin{equ}
F_\eps(t) := \int_0^t \scal{G_\eps(s),f(s)}\,ds\;.
\end{equ}
Assume furthermore that, for every $\kappa>0$ and every $p > 0$, there exists a constant $C$ such that
\begin{equ}[e:boundsup]
\E \sup_{t \in [0,1]} \|G_\eps(t)\|^p \le C \eps^{-\kappa}\;.
\end{equ}
Then, for every $\gamma < 2\tilde\alpha / (1 + 2\tilde\alpha)$, there exists a constant $C$ depending
only on $p$ and $\gamma$  such that 
\begin{equ}
\E \|F_\eps\|_{C^{1-\tilde\alpha}}^p \le C \bigl(\E \|f\|_{C^{\tilde\alpha}}^{2p}\bigr)^{1/2}
\eps^{\gamma p}\;,
\end{equ}
where $\|\cdot\|_{C^{\tilde\alpha}}$ denotes the $\tilde\alpha$-H\"older norm for $\KK$-valued
functions on $[0,1]$.
\end{proposition}

Note that, if we can choose $\tilde\alpha<1/2$, but arbitrarily close, then
we can choose $\gamma<1/2$, but arbitrarily close, too.

\begin{proof}
We focus only on the H\"older part of the norm. The $L^\infty$ part follows easily, as 
$F_\eps(0)=0$, e.g. by taking $s=0$ in the following.

Choose $\delta > 0$ to be fixed later.  Moreover, for every pair $s$ and $t$ in $[0,1]$, set $\bar s =
\delta[\delta^{-1}s]$ and $\bar t = \delta[\delta^{-1}t]$. We furthermore define
$f_\delta(t) = f(\bar t)$. One then has
\begin{equs}
\Bigl|\int_s^t \scal{G_\eps(r),f(r)}\,dr\Bigr| &\le \Bigl|\int_s^t
\scal{G_\eps(r),\bigl(f(r)-f_\delta(r)\bigr)}\,dr\Bigr|
+\Bigl|\int_s^t \scal{G_\eps(r), f_\delta (r)}\,dr\Bigr| \\
&\le \|G_\eps\|_{L^\infty(\KK)} \delta^{\tilde\alpha} |t-s| \|f\|_{C^{\tilde\alpha}} + \Bigl|\int_s^t \scal{G_\eps(r),
f_\delta (r)}\,dr\Bigr|\;.
\end{equs}
The second term in the right hand side can be bounded by
\begin{equ}
\One_{|t-s| \ge \delta}\,\Bigl|\int_{\bar s}^{\bar t}\scal{G_\eps(r), f_\delta (r)}\,dr\Bigr| +
\min\{|t-s|,2\delta\} \|f\|_{C^{\tilde\alpha}} \|G_\eps\|_{L^\infty(\KK)}\;.
\end{equ}
The first term of this expression is in turn bounded by
\begin{equ}
{2|t-s| \over \delta}\Bigl(\sup_{n < \delta^{-1}}\Bigl\| \int_{n\delta}^{(n+1)\delta}
G_\eps(r) \,dr\Bigr\|\Bigr) \|f\|_{C^{\tilde\alpha}}\;.
\end{equ}
Collecting all these expressions yields
\begin{equ}
\|F_\eps\|_{C^{1-\tilde\alpha}}  \le C \|G_\eps\|_{L^\infty(\KK)}\|f\|_{C^{\tilde\alpha}}  \delta^{\tilde\alpha} +
\delta^{-1}\Bigl(\sup_{n < \delta^{-1}}\Bigl\| \int_{n\delta}^{(n+1)\delta}
G_\eps(r)\,dr\Bigr\|\Bigr) \|f\|_{C^{\tilde\alpha}}\;.
\end{equ}
Choosing $\delta = \eps^{2/(1+2\tilde\alpha)}$, applying Lemma~\ref{lem:supav}, and using \eref{e:boundsup}
easily concludes the proof.
\end{proof}

We are actually going to use the following corollary of Proposition~\ref{prop:homog}:
\begin{corollary}\label{cor:homog}
Let $\hat z$ be as above and let $\alpha$ be as in Assumptions~\ref{ass:3} and \ref{ass:4}.
Fix $T>0$ and let $f_i$ with $i \in \{1,2,3\}$ be $\tilde \alpha$-H\"older continuous functions on $[0,T]$
with values in
$\bigl((\HH^\alpha)^{\otimes i}\bigr)^*$, respectively.
Let $F_\eps$ be given by
\begin{equ}
F_\eps(t) := \int_0^t \Bigl(\bigl(f_1(s)\bigr)(\hat z) + \bigl(f_2(s)\bigr)(\hat z\otimes \hat z - \widehat Q)
+ \bigl(f_3(s)\bigr) (\hat z\otimes \hat z\otimes \hat z) \Bigr) \,ds\;.
\end{equ}
Then, for every $\gamma < 2\tilde\alpha / (1 + 2\tilde\alpha)$, there exists a constant $C$ depending
only on $p$ and $\gamma$  such that 
\begin{equ}
\E \sup_{t \in [0,T]} |F_\eps(t)|^p \le C \eps^{\gamma p} \bigl(\E \bigl(\|f_1\|_{C^{\tilde\alpha}} + \|f_2\|_{C^{\tilde\alpha}} + \|f_3\|_{C^{\tilde\alpha}}\bigr)^{2p}\bigr)^{1/2}
\;,
\end{equ}
where $\|\cdot\|_{C^{\tilde\alpha}}$ denotes the $\tilde\alpha$-H\"older norm for $\bigl((\HH^\alpha)^{\otimes i}\bigr)^*$-valued 
functions on $[0,\tau^*]$.
\end{corollary}

\begin{proof}
Note that $\hat z$ satisfies \eref{e:boundsup} with $\KK =\CH^\alpha$.
This follows for example from the proof of \cite[Theorem~5.9]{DapZab92}.
The statement is then a consequence of Corollary~\ref{cor:averagez}
and of Proposition~\ref{prop:homog}.
\end{proof}

%
\subsection{The reduction of the slow modes}
%
%
%
We now use the results of the previous subsection in order to show that most of the terms
that appear on the right hand side of eqn.~\eqref{e:ampl1} are of order $\cO(\eps^{1/2-})$. Note 
first that we can replace  all occurrences of $z$ in \eref{e:ampl1} by the stationary 
process  $\hat z$ without changing the order of magnitude of the remaining term $R$. 
We are now ready to state the main result of this section.
\begin{proposition}
\label{prop:sndstep} 
Under Assumptions~\ref{ass:1}--\ref{ass:4} and with $x$ and $\hat z$ defined as above,
we obtain
\begin{equs}
x(t) &= x(0)
+ \int_0^t  P_c B (I\otimes_s L)^{-1}\bigl(\hat z(s) \otimes_sQ\, dW(s)\bigr)\\
&\quad + 2\int_0^t P_c B\bigl(x(s), L^{-1} Q\,dW(s)\bigr) \\
&\quad + \tilde \nu  \int_0^t x(s)\,ds - \tilde{\eta} \int_0^t \scal{e_1,x(s)}^3 e_1\,ds 
+ R(t)\; , 
\label{e:prop-2nd}
\end{equs}
where  $\|R(\cdot)\|=\cO(\eps^{1/2-})$ and the constants $\tilde \nu$ and $\tilde \eta$ 
are defined as
\begin{equs}
\tilde\nu &= \nu + 2\scal[b]{e_1,B \bigl( (I\otimes_s L)^{-1} (B_s \otimes_s I)
+  (I \otimes L^{-1}B_s)
\\ &\quad +  2(B_c \otimes L^{-1})\bigr)(e_1 \otimes \widehat Q)} \;, 
\label{e:nu} \\
\tilde \eta &= -2\scal{e_1,B(e_1,L^{-1} B_s(e_1,e_1))} \;. \label{e:eta}
\end{equs}
Here, we used the notation $B_s := P_s B$ and $B_c := P_c B$.
\end{proposition}

\begin{proof}
First we replace all instances of $z$ by $\hat z$ on the
right hand side of eqn.~\eqref{e:ampl1}, which results in an error of order $\CO(\eps^{1-})$ which is absorbed into $R$.
This is a straightforward but somewhat lengthy calculation which we do not reproduce here.
We rely on
$$
z(t)-\hat{z}(t) = e^{-tL\eps^{-2}} (P_s v_0- z(0)) = e^{-tL\eps^{-2}} \cO(\eps^{0-})\;. 
$$
Note that obviously $\|z-\hat{z}\|\not=\cO(\eps^{1-})$, as bounds uniformly in time are not available
due to transient effects on timescales smaller than $\cO(\eps^2)$.
Nevertheless, we only bound the error in integrated form, 
which is sufficient for our application. 
Actually, it is possible to check that the error terms which result from this 
substitution are of $\cO(\eps^{2-})$. The only exception to this is the 
stochastic integral, where we apply the Burkholder-Davies-Gundy inequality in order to
get a remainder term of $\cO(\eps^{1-})$.

Once this substitution has been performed, the proposition follows from an application 
of Corollary \ref{cor:homog} to the modified eqn.~\eqref{e:ampl1}.
 The fact that, for every $\tilde \alpha > 1/2$, the various integrands are indeed
$\tilde \alpha$-H\"older continuous functions with values in
$\bigl((\HH^\alpha)^{\otimes i}\bigr)^*$ and H\"older norm of order $\CO(\eps^{0-})$
follows from Corollary~\ref{rem:Hoelder}.
\end{proof}

%
%
\section{Approximation of the Martingale Term}
\label{sec:fin-red}
%
%

This section deals with the final reduction step for  the general system \eref{e:x-v}. We start 
by eliminating the stochastic integral of the type $\int_0^t \hat{z} \otimes Q dW(s)$ from 
\eqref{e:prop-2nd}. In fact, we show that we can replace the martingale part
in eqn.~\eqref{e:prop-2nd} by a single stochastic integral of the type
\begin{equ}
\int_0^t \sqrt{\sigma_a + \sigma_b a^2(s)}\,dB(s)\;,
\end{equ}
against a \emph{one-dimensional} Wiener process $B$. Note that this section is superfluous in the
particular case where the first stochastic integral in \eref{e:prop-2nd} vanishes. This is
the case for example in the situation considered in \cite{Ro:03}. See Theorem \ref{thm:burgers1} in
the next section.

We emphasise that although all the previous steps 
are easily extended to higher dimensions, this step is valid only under the assumption 
that the kernel of $L$ is one-dimensional.
 
%
\subsection{An abstract martingale approximation result}
\label{susec:approx}
%
We start with the following lemma; we will use it to approximate the martingale part of 
equation \eref{e:prop-2nd} by a stochastic integral against a one dimensional Brownian
motion.

\begin{lemma}\label{lem:mart}
Let $M(t)$ be a continuous $\cF_t$-martingale with quadratic variation $f$ and let $g$ be 
an arbitrary $\cF_t$-adapted increasing process with $g(0) = 0$. Then, there exists
a filtration $\tilde \cF_t$ with $\cF_t \subset \tilde\cF_t$ and a continuous 
$\tilde \cF_t$-martingale $\tilde M(t)$ with quadratic variation $g$ such that, for every 
$\gamma < 1/2$ there exists a constant $C$ with
\begin{equs}
\E \sup_{t \in [0,T]} |M(t) - \tilde M(t)|^p &\le C \bigl(\E g(T)^{2p}\bigr)^{1/4}
\bigl(\E \sup_{t \in [0,T]} |f(t) - g(t)|^{p}\bigr)^{\gamma}\\
& \qquad + C \E \sup_{t \in [0,T]} |f(t) - g(t)|^{p/2} \;.
\end{equs}
\end{lemma}

\begin{proof}
Define the adapted increasing process $h$ by
\begin{equ}
h(t) = f(t) + \inf_{s \le t} \bigl(g(s) - f(s)\bigr)\;.
\end{equ}
note that one has $h \le g$ almost surely. Furthermore, one has
$$
0\le h(t) - h(s) \le f(t) - f(s)
$$ 
for every $t \ge s$, so that one has $0¸\le {dh \over df} \le 1$ almost
surely. Define a martingale $\hat{M}(t)$ with quadratic variation $h$ by the It\^o integral
\begin{equ}
 \hat{M}(t) = \int_0^t \sqrt{{dh \over df}(s)}\,dM(s)\;.
\end{equ}
Define now an increasing sequence of random times $T_t$ by
\begin{equ}
T_t =\inf \bigl\{s\ge 0\,|\, h(s) \ge g(t)\bigr\} \ge t \;.
\end{equ}
Note that since $h \le g$ almost surely, the times $T_t$ are actually stopping times with respect to $\cF_t$,
so that the time-changed process $\tilde M(t) = \hat M(T_t)$ is a martingale with
quadratic variation $g$. Note that $\tilde M(t)$ is a martingale with respect to the filtration $\tilde \cF_t$ induced 
by the stopping times $T_t$. Note also that $\cF_t \subset \tilde \cF_t$ as a consequence of the fact that $T_t \ge t$
almost surely.

It remains to show that $\tilde M$ satisfies the required bound. Let us
start by defining the martingale $\Delta$ as the difference $\Delta = M - \hat M$.
The quadratic variation $\scal{\Delta}$ of $\Delta$ is then bounded by
\begin{eqnarray*}
\scal{\Delta}(t) 
&=& \int_0^t \left(1 - \sqrt{{dh \over df}(s)} \right)^2\,df(s) 
\\ &&\le \int_0^t \Bigl(1 - {dh \over df}(s)\Bigr)\,df(s) 
= f(t) - h(t)\\
&=& \sup_{s \le t} \bigl(f(s) - g(s)\bigr)\;.
\end{eqnarray*}
Applying the Burkholder-Davies-Gundy
inequalities {\cite[Cor. IV.(4.2)]{RevYor99}} 
to this bound yields the existence
of a universal constant $C$ such that
\begin{equ}
\E \sup_{t \in [0,T]} |\Delta(t)|^p \le C \E \sup_{t \in [0,T]}|f(t) - g(t)|^{p/2}\;.
\end{equ}

Before we turn to bounding the difference between $\hat M$ and $\tilde M$, we
show that if $F$ is an arbitrary positive random variable, $B$ is a Brownian
motion, $\gamma < 1/2$, and $q > p > 1$, then there exists a constant $C$ depending only on $p$, $q$ and $\gamma$,
 such that
\begin{equ}[e:boundHolder]
\E \|B\|_{\gamma,F}^p
\le C \bigl(\E F^q\bigr)^{p/2q}\;,
\end{equ}
where we defined
\begin{equ}
\|B\|_{\gamma,F} = \sup_{0 \le s < t \le F} {|B(t) - B(s)|\over |t-s|^{\gamma}}\;.
\end{equ}
One has indeed for every $K>0$ and every $L > 0$ the bound
\begin{equ}
\P \bigl(\|B\|_{\gamma,F} > K\bigr)
\le \P(F \ge L) + \P \bigl(\|B\|_{\gamma,L} > K\bigr)\;.
\end{equ}
Applying Chebyshev's inequality and using the Brownian scaling together with the fact that
the $\gamma$-H\"older norm of a Brownian motion on $[0,1]$ has moments of all orders, this yields
for every $q > 0$ the existence of a constant $C$ such that
\begin{equ}
\P \bigl(\|B\|_{\gamma,F} > K\bigr)
\le \inf_{L>0} \Bigl({\E F^q \over L^q} + {\E \|B\|_{\gamma,1}^{2q} L^{q} \over K^{2q}}\Bigr)
\le C{(\E F^q)^{1/2} \over K^q}
\;.
\end{equ}
The bound \eref{e:boundHolder} then follows immediately from the fact that if
a positive random variable $X$ satisfies $\P(X>K) \le (a/K)^q$ for some $a$, some
$q$ and every $K>0$ then, for every $p < q$,  there exists a constant $C$ such
that $\E |X|^p \le C a^p$.

Note now that it follows from our construction that
there exists a Brownian motion $B$ such that
$\hat M(t) = B(h(t))$ and $\tilde M(t) = B(g(t))$. Noting that $h \le g$ and setting $G = g(T)$,
we have
\begin{equs}
\E \sup_{t \in [0,T]} |\hat M(t) - \tilde M(t)|^p 
&\le \E \|B\|_{\gamma,G}^p \sup_{t \in [0,T]}|h(t) - g(t)|^{\gamma p}
\\
&\le \E \|B\|_{\gamma,G}^p \sup_{t \in [0,T]} |f(t) - g(t)|^{\gamma p}
\end{equs}
and the result follows from \eref{e:boundHolder} and Young's inequality.
\end{proof}
%
%
%
\subsection{Application to the SPDE}
\label{sec:secred2}
%
%

Before we state the next result, we introduce some notation. Let $\gamma \in \HH$ and
$\Gamma \colon \HH^{\alpha} \to \HH$ be defined by
\begin{equ}
\scal{y,\gamma} = 2\scal[b]{e_1, B(e_1, L^{-1}Q y)}\;,\quad
\scal{y,\Gamma z} = \scal[b]{e_1, B (I \otimes_s L)^{-1} \bigl(z \otimes Qy\bigr)}\;.
\end{equ}
The facts that $\gamma \in \HH$ and $\Gamma$ is bounded follow  from
Lemma \ref{lem:3.8} together with the fact that
Assumption~\ref{ass:4} implies in particular that $Q$ is a bounded operator from $\HH$ to $\HH^{\alpha-1}$.
Note that $\Gamma$ is actually bounded as an operator from $\HH^{\alpha-1}$ to $\HH$, but we
will not need this fact.

\begin{theorem}\label{theo:main}
Let Assumptions \ref{ass:1}--\ref{ass:3} hold
and let $(x(t),y(t))$ be the solution of \eqref{e:x-v}. Let $\tilde \nu, \, \tilde \eta$ 
be given by \eqref{e:nu} and \eqref{e:eta}, respectively, and define
\begin{equ}[e:ab]
\sigma_a = \|\gamma\|^2\;,\qquad \sigma_b = \tr \bigl(\Gamma \widehat Q \Gamma^*\bigr)\;, 
\end{equ}
where we identify $\widehat Q$ from \eref{e:defQhat}
 with the corresponding operator\footnote{An element $u\otimes v$ of $\HH^\alpha\otimes\HH^\alpha$
defines an operator from $(\HH^\alpha)^*$ to $\HH^\alpha$ by $(u\otimes v)(f)= u \scal{f,v}$} 
from $(\HH^\alpha)^*$ to $\HH^\alpha$. 

Define finally $X(t) = \langle x(t), e_1 \rangle$. Then, there exists a Brownian motion 
$B$ such that, if $a$ is the solution to the SDE
\begin{equ}[e:ampl]
da(t) = \tilde \nu a(t) - \tilde\eta a^3(t) + \sqrt{\sigma_b  + \sigma_a a^2(t)}\,dB(t), 
\quad a(0) = X(0),
\end{equ}
then, for every $p > 0$ and every $\kappa > 0$, there exists a constant $C$ such that
\begin{equation*}
\E \sup_{t \in [0,\tau^*]} |X(t) - a(t)|^p \le C \eps^{p/4 - \kappa}\;,
\end{equation*}
for every $\eps \in(0,1)$, where $\tau^*$ is defined in \eref{tau*-def}.
\end{theorem}
\begin{proof} From 
Proposition \ref{prop:sndstep} we have that, with the notations introduced above,
\begin{eqnarray*}
X(t) &=& X(0) +\tilde \nu  \int_0^t X(s)\,ds - \tilde{\eta} \int_0^t X^3(s)\,ds 
\\ && + \int_0^t  \scal[b]{\Gamma\hat z(s), dW(s)\bigr)}
+ \int_0^t X(s) \scal[b]{\gamma,dW(s)} +  R_2(t)\; , 
\end{eqnarray*}
where $R_2=\cO(\eps^{1/2- \kappa})$. Denote by $M(t)$ the martingale
\begin{equ}
M(t) =\int_0^t  \scal[b]{\Gamma \hat z(s), dW(s)\bigr)}
+ \int_0^t X(s) \scal[b]{\gamma,dW(s)} \;.
\end{equ}
Its quadratic variation is given by
\begin{equ}[e:qvf]
f(t) = \int_0^t \|\gamma  X(s) + \Gamma \hat z(s)\|^2\,ds\;. 
\end{equ}
It now follows from Corollary~\ref{cor:homog} that
\begin{equ}[e:qvg]
|f(\cdot) - g(\cdot)| = \CO(\eps^{1/2 -})\;,
\quad\text{where}\quad
 g(t) = \int_0^t \bigl(\sigma_a X^2(s) + \sigma_b \bigr)\,ds\;.
\end{equ}

Denote by $\tilde{M}(t)$ the martingale with quadratic variation $g(t)$ given by 
Lemma~\ref{lem:mart}
and by $\tilde x$ the solution to
\begin{equ}
d\tilde x = \tilde \nu \tilde x \,dt - \tilde \eta \tilde x^3\,dt + d\tilde M(t)\;,
\qquad \tilde x(0) = x(0)\;.
\end{equ}
It follows from Lemma~\ref{lem:mart} that
$M(t) - \tilde M(t) = \CO(\eps^{1/4 -})$.
Therefore, using a standard estimate stated below in Lemma \ref{lem:a-priori},
\begin{equ}
\label{e:x-x}
x(t) - \tilde x(t) = \CO(\eps^{1/4 -})\;.
\end{equ}
The martingale representation theorem \cite[Thm V.3.9]{RevYor99} ensures that one can
enlarge the original probability space in such a way that  there exists
a filtration $\tilde \cF_t$, and an $\tilde \cF_t$-Brownian motion
$B(t)$,
such that both $x(t)$ and $\tilde{x}(t)$ are $\tilde \cF_t$-adapted and such that
\begin{equ}
d \tilde{x}(t) = \tilde\nu\tilde{x}(t) \,dt - \tilde\eta \tilde{x}^3(t)\,dt + 
\sqrt{\sigma_b + \sigma_a X^2(t)}\,d B(t)\;.
\end{equ}
Note that  in general the $\sigma$-algebra $\tilde \cF_t$ is strictly larger than the one
generated by the Wiener process $W$ up to time $t$. This is a consequence of the construction
of Lemma~\ref{lem:mart} where one has to `look into the future' in order to
construct $\tilde M$.

We finally define the process $a$ as the solution to the SDE
\begin{equ}
da(t) = \tilde\nu a(t) \,dt - \tilde\eta a^3(t)\,dt + \sqrt{\sigma_b + \sigma_a  a^2(t)}
\,d B(t)\;.
\end{equ}
Denote $\rho = a - \tilde x$ and $G(x) = \sqrt{\sigma_b + \sigma_a x^2}$. Then, one has
\begin{equ}
d\rho^2(t) \le 2\tilde \nu \rho^2(t)\,dt
 + \bigl(G(a(t)) - G(X(t))\bigr)^2\,dt+ 2\rho \bigl(G(a(t)) - G(X(t))\bigr)\,d B(t)\;.
\end{equ}
Using the fact that $G$ is globally Lipschitz, this yields the existence of a
constant $C$ such that
\begin{equ}
d\rho^2(t) \le C \rho^2(t)\,dt
 + C |X(t) - \tilde{x}(t)|^2\,dt+ 2\rho \bigl(G(a(t)) - G(X(t))\bigr)\,d B(t)\;.
\end{equ}
It is now easy to verify, using It\^o's formula, \eref{e:x-x}, and 
the Burkholder-Davies-Gundy inequality, that
\begin{equ}
\rho=\CO(\eps^{1/4 -})
\quad\text{and thus}\quad
X(t) - a(t)=\CO(\eps^{1/4 -})\;,
\end{equ}
which is the required result.
\end{proof}

Let us finally state the a-priori estimate used in the previous proof.

\begin{lemma}
\label{lem:a-priori}
Fix $\nu\in\R$ and $\eta>0$.
Let $M_i(t)$ be martingales (not necessary with respect to the same filtration),
and $x_i(t), \, i=1,2$ be solution of the following SDEs
\begin{equation}\label{e:xi}
d x_i(t) = \nu x_i(t) \,dt - \eta x^3_i(t)\,dt + d M_i(t)\;,
\end{equation}
with $x_1(0)-x_2(0)=\cO(\eps^{1/4-})$ and  $x_1(0)=\cO(\eps^{0-})$.
Suppose furthermore $M_i(t)=\cO(\eps^{0-})$ and $M_1(t)-M_2(t)=\cO(\eps^{1/4-})$,
then 
$$x_1(t)-x_2(t)=\cO(\eps^{1/4-})\;.$$
\end{lemma}

\proof This is a straightforward a priori estimate which relies on the stable cubic
nonlinearity in \eqref{e:xi}. First, one easily sees from It\^o's formula that 
$x_i(t)=\cO(\eps^{0-})$. Then using the transformation $\hat x_i(t)=x_i(t)-M_i(t)$ to 
random ODEs for $\hat x_i(t)$, we can write down an ODE for the difference 
$\hat x_1(t)-\hat x_2(t)$, which we can bound pathwise by direct a priori estimates. 
We will omit the details. 
\qed

%
\subsection{Main Result}
%

Let us finally put the results obtained in this and the previous two sections together to
obtain our final result for the system of SDEs \eqref{e:x-v}. 
\begin{theorem}
\label{thm:final}
Let Assumptions \ref{ass:1}--\ref{ass:4} be true.
Let $(x(t), \, y(t))$ be a solution of \eref{e:x-v}. Furthermore, let $z(t)$ be the 
OU-process defined in \eref{e:ou} and $\tau^*$ the stopping time from 
\eref{tau*-def}. Let finally $\sigma_a$, $\sigma_b$,  $\tilde \nu$, and 
$\tilde\eta$ be defined in \eref{e:ab}, \eref{e:nu}, and \eref{e:eta} 
respectively. Then there exists a Brownian motion $B(t)$ such that, if $a(t)$
is a solution of 
\begin{equ}[e:ampla]
d a(t) = \tilde \nu a(t) - \tilde\eta a^3(t) + \sqrt{\sigma_b + \sigma_a a^2(t)}\,dB(t)\;, 
\quad a(0)=\langle x(0),e_1 \rangle \;, 
\end{equ}
then for all $T>0$, $R>0$, $p>0$ and $\kappa>0$ there is a constant $C$ such that
for all $\eps\in(0,1)$ and all $\|x(0)\|_\alpha <R$ and $\|y(0)\|_\alpha<R$ we have that
$$
\P (\tau^*>T) > 1-C\eps^p\;,
$$
$$
\E \sup_{t\in[0,\tau^*]} \|x(t) - a(t)e_1 \|^p_\alpha \le C\eps^{p/4-\kappa}\;,
\quad\text{and}\quad
\E \sup_{t\in[0,\tau^*]}\|y(t) - z(t)\|^p_\alpha \le C \eps^{p-\kappa}\;.
$$
\end{theorem}
\begin{proof}
The approximation of $y(t)$ by $z(t)$ is already verified in Corollary \ref{cor:1stBtau*} 
and Lemma \ref{lem:approxOU}. The approximation of $x(t)$ by $a(t)$ follows from Theorem 
\ref{theo:main}. The bound on the stopping time $\tau^*$  follows then easily from the 
fact that $a(t)=\CO(1)$ and $z(t)=\CO(1)$ uniformly on $[0,T]$.
\end{proof}

\begin{remark}\label{rem:scaling}
With the notation $B_{k \ell m} = \scal{B(e_k, e_{\ell}), e_m }$,
formulas \eref{e:ab}, \eref{e:nu}, and \eref{e:eta} for the coefficients in the
amplitude equation can be written in the form
\minilab{e:coeffs}
\begin{equs}
\tilde \nu &= \nu + \sum_{k=2}^\infty  {2 B_{k11}^2 q_k^2 \over \lambda_k^2}
    +  \sum_{k,\ell=2}^\infty\frac{B_{k11} B_{\ell \ell k} q_\ell^2}{\lambda_k \lambda_\ell}
	+  \sum_{k,\ell=2}^\infty 
		\frac{2B_{k\ell 1}B_{k1\ell}}{ \lambda_k+\lambda_\ell}
		\frac{q_k^2}{ \lambda_k} \;, \qquad\label{e:coeff1}\\
\tilde \eta &= -\sum_{k=2}^{\infty} \frac{2 B_{k11} B_{11k}}{\lambda_k}\;,\label{e:coeff2}\\
\sigma_a &= \sum_{k=2}^\infty\frac{4 B_{k11}^2 q_k^2}{\lambda_k^2}\;,\quad
\sigma_b = \sum_{m,k=2}^{\infty} \frac{2 B_{km1}^2 q_k^2 q_m^2}{(\lambda_k+\lambda_m)^2\lambda_k}
\;.\label{e:coeff3}
\end{equs}
If one chooses to expand the solution in a basis which is not normalised, \ie one takes
basis vectors $\tilde e_k = c_k e_k$, then the coefficients appearing in the right-hand side
of the equation for the expansion transform according to
\begin{equ}
\tilde B_{k\ell m} = B_{k\ell m} {c_k c_\ell\over c_m} \;,\qquad \tilde q_k = {q_k \over c_k}\;.
\end{equ}
It is straightforward to see that $\tilde \nu$ and $\sigma_a$ are unchanged under this transformation,
whereas $\tilde \eta$ is mapped to $c_1^2 \tilde \eta$ and $\sigma_b$ is mapped to $\sigma_b / c_1^2$
as expected.
\end{remark}

%
%
\section{Application: The Stochastic Burgers Equation}
\label{sec:burgers}
%
%
In this section we apply our results to a modified stochastic Burgers equation:
\begin{equation}\label{e:burgers_appl} 
d u = (\partial_x^2 +1) u \, dt + u \partial_x u \, dt + \eps^2 \nu u \, dt + \eps \, Q \,dW,
\end{equation}
on the interval $[0, \pi]$, with Dirichlet boundary conditions. 
We take
\begin{equs}
\HH =L^2([0,\pi])\;, \quad e_k(x)=\sqrt{2\over\pi}\sin(kx)\;, \qquad \\ 
B(u,v)=\frac12\partial_x (uv)\;, \quad L=-\partial_x^2-1\;,\quad \lambda_k=k^2-1\;.
\end{equs}
We also take $W(t)$ to be a cylindrical Wiener process on $\HH$
and $Q$ a bounded operator with $Qe_1 = 0$ and $Q e_k = q_k e_k$ for $k \ge 2$.
It follows that $\HH^\alpha=H^\alpha_0([0,\pi])$ is the standard fractional Sobolev space
defined by the Dirichlet Laplacian on $[0,\pi]$.

With this choice, using the notation from Remark~\ref{rem:scaling}, we get
\begin{equation}
\label{e:repBklm}
B_{k \ell m} 
= \frac{1}{2\sqrt{2\pi}} \big( |k+ \ell| \delta_{k + \ell, m} - |k - \ell| \delta_{|k -\ell|,m} \big)\;,
\end{equation}
where $\delta_{k\ell}$ is the Kronecker delta symbol.

It is possible to check that assumptions \ref{ass:1}, \ref{ass:2}, and  
\ref{ass:3} are satisfied for any $\alpha \ge 0$ since, for smooth functions
$u$, $v$, and $w$, one has for example
\begin{equs}
\Bigl|\int_0^\pi (uv)'(x)\, w(x)\, dx\Bigr|
&= \Bigl|\int_0^\pi u(x)v(x)w'(x)\, dx\Bigr|
\le \|u\| \|v\| \|w'\|_{\L^\infty}\\
& \le C\|u\| \|v\| \|w\|_{\HH^{\gamma}}\;,
\end{equs}
provided that one takes $\gamma < -3/2$. (Values of $\alpha$ other than $0$ can be obtained in
a similar way by using different Sobolev embeddings, see also \cite{DP-D-T:94}.)
Whether the trace-class assumption on $Q^2 L^{\alpha-1}$ is
satisfied or not depends of course in a crucial way on the coefficients $\{q_k \}_{k=1}^{\infty}$. 

The following result justifies the formal asymptotic calculations presented in \cite{Ro:03}.
\begin{theorem}\label{thm:burgers1}
Let $u$ be a continuous $H^1_0([0,\pi])$-valued solution  of \eref{e:burgers_appl}
with initial condition $u(0)=\cO(\eps)$,
and assume that the driving noise $W$ is given by $\sigma \sin(2x) w(t)$ for a standard
one-dimensional Wiener process $w$. 
Then there are Brownian motions $B(t)$ and 
$\beta(t)$ (not necessarily adapted to the same filtration) such that
if $a$ is the solution of
$$ da=\Bigl(\nu-\frac{\sigma^2}{88}\Bigr)\,a\; dt-\frac1{12} a^3 \; dt + \frac{\sigma}{6} 
|a|\circ dB, \quad \eps a(0)=\frac{2}{\pi} \big( u(0), \sin(\cdot) \big)_{L^2}
$$
and
$$ R(t) =  \frac{1}{\eps}
e^{-Lt} P_s u(0) + \left( \int_0^t e^{-3(t-s)}d\beta(s) \right) \sin(2 \cdot),
$$
then for all $\kappa,p>0$ there is a constant $C$ such that
$$
\E \Bigl(\sup_{t\in[0,T\eps^{-2}]}\|u(t)-\eps a(\eps^2t)\sin(\cdot) -\eps R(t) 
 \|_{H^1}^p \Bigr) \leq C \eps^{\frac{3p}{2}-\kappa}  
$$
\end{theorem}
\begin{proof}
Note first that Assumption \ref{ass:4} is true for all $\alpha$ is this case, so that
all the assumptions of Theorem~\ref{thm:final} are satisfied. Furthermore, we
can use formulas \eqref{e:coeffs} to obtain
\footnote{Notice that $a$ is the amplitude of the mode $\sin(x)$ which is not normalized. 
This is in order to be consistent with earlier works on the stochastic Burgers equation. 
The modification of the formulas for the constants that 
appear in the amplitude equation in this situation is given by Remark~\ref{rem:scaling}.}, 
$$
\tilde{\eta} =  \frac{1}{12}
\;,\quad
\sigma_a =  \frac{\sigma^2}{36}
\;,\quad
\sigma_b=0\;,\quad
\tilde{\nu}=\nu+\frac{\sigma^2}{72}-\frac{\sigma^2}{88}\;.
$$
Note that the second term in the expression for $\tilde\nu$ gives the It\^o-Stratonovich correction.

However, the claim does not follow immediately, 
since we wish to get an error estimate of order $\eps^{3/2}$ instead of $\eps^{5/4}$.
Retracing the proof of Theorem~\ref{thm:final}, we see that the claim follows if we can
show that $|f-g| = \CO(\eps^{-})$, where $f$ and $g$ are as in \eref{e:qvf} and \eref{e:qvg}. 
In our particular case, one has $\gamma = 0$, so that
\begin{equ}
f(t) = \int_0^t \|\Gamma \hat z(s)\|^2\,ds\;,\qquad 
 g(t) =  \sigma_b t\;.
\end{equ}
The result now follows from Lemma~\ref{lem:qvsimple} below.
\end{proof}

\begin{lemma}\label{lem:qvsimple}
Let $\hat z$ be as in Corollary~\ref{cor:averagez}. Then, for every final time $T$, every $p>0$ and every $\kappa > 0$ there
exists a constant $C$ such that
\begin{equ}[e:boundqv]
\E \sup_{t \in [0,T]} \Bigl\|\int_0^t \bigl(\hat z(r)\otimes \hat z(r)- \widehat Q\bigr) \,dr\Bigr\|_\alpha^{2p} \le C \eps^{2p - \kappa}\;.
\end{equ}
\end{lemma}

\begin{proof}
We subdivide the interval $[0,T]$ into $N$ subintervals of length $T/N$, and we use the notation
$t_k = k T/N$.
Using exactly the same argument as in the proof of Lemma~\ref{lem:supav}, we see that, for every
$p>0$ and every $\kappa > 0$ there exists a constant $C$ such that
\begin{equ}
\E \sup_{k \in \{0,\ldots,N\}} \Bigl\|\int_0^{t_k} \bigl(\hat z(r)\otimes \hat z(r)- \widehat Q\bigr) \,dr\Bigr\|_\alpha^{2p} \le C N^\kappa \eps^{2p}\;.
\end{equ}
On the other hand, we know from Lemma~\ref{lem:bound-y} that the $\HH^\alpha$-norm of
the integrand in \eref{e:boundqv} is of order $\CO(\eps^{0-})$ uniformly in time. The claim
then follows by taking $N \approx \eps^{-1}$.
\end{proof}
%
%
\begin{remark}
In the case where only the second mode is forced by noise, one can actually take 
$\beta(t) = B(t)$, and $\beta(t)$ could be chosen to be a rescaled version of the Brownian motion that 
appears in equation \eqref{e:burgers_appl}. 
\end{remark}

The following theorem covers the case where $W(t)$ is space-time white noise which is 
constrained to be antisymmetric around $x = \pi$.
This corresponds to the case $q_k=\sigma$ for all $k\ge2$.
It is easy to check that Assumption \ref{ass:4} is satisfied for all $\alpha<\frac12$.

We again use Remark~\ref{rem:scaling} with $c_k = \sqrt{\pi / 2}$ in order to
compute the coefficients. We obtain
\begin{equs}
\tilde{\eta} &= \frac{1}{12}
\;,\quad
\sigma_a =  \frac{\sigma^2}{18\pi}\;,\\
\sigma_b &= c_b \sigma^4\;,\quad
c_b = 
\frac{1}{2\pi^2}
 \sum_{k=2}^\infty \frac{1}{(2k^2 + 2k + 1)(k^2-1)(k^2 + 2k)}\;,
\end{equs}
and finally, 
\begin{equs}
\tilde{\nu}-\nu &= \frac{\sigma^2}{36\pi} - \frac{\sigma^2}{4\pi} 
 \sum_{k=2}^\infty \Big(\frac{1}{k-1}  
- \frac{1}{k(k+1)}   \Big)\frac{1}{2k^2+2k-1} \;.
\end{equs}
Note again that the first term in this expression is the Stratonovitch correction for the
multiplicative noise term. This finally leads to

\begin{theorem}
Assume that $\alpha\in[0,\frac12)$ and let $u$ be a continuous 
$H^\alpha_0([0,\pi])$-valued solution of \eref{e:burgers_appl} with initial condition 
$u(0)=\cO(\eps)$. Assume furthermore that the covariance of the noise satisfies  
$q_k = \sigma$  for $k \geq 2$. Then there is a Brownian motion $B(t)$  
(not necessarily adapted to the filtration of $W(t)$) such that if $a(t)$ is a solution of 
$$d a(t) = \tilde \nu a(t) - \tilde\eta a^3(t) + \sqrt{\sigma_b + \sigma_a a^2(t)}\,dB(t), \quad \eps a(0)=\frac{2}{\pi} \big(u(0),\sin(\cdot) \big)_{L^2}$$
where the constants are defined above, and 
$$  R(t) =   \frac{1}{\eps}e^{-tL} P_s u(0) + \int_0^t e^{-(t-s)L} QdW(s)\;,
$$
then for all $\kappa,p>0$ there is a constant $C$ such that
$$
\E \Bigl( \sup_{t\in[0,T\eps^{-2}]}\|u(t)-\eps a(\eps^2t)\sin(\cdot) -\eps R(t)\|^p_{H_0^\alpha}
\Bigr)  \leq  C \eps^{\frac{5p}{4}-\kappa} 
$$
for $\eps$ sufficiently small.
\end{theorem}
%

\bibliography{mybib,newrefs}
\bibliographystyle{Martin}
\end{document}